\documentclass[12pt,a4paper]{article}
\usepackage{amsmath,amssymb,amsthm,array,amsfonts}
\usepackage{mathrsfs} 

\fussy
\begin{document}
\fontshape{r}
\selectfont 
\begin{center}
\textbf{ MINIMAX STATE ESTIMATION FOR A DYNAMIC SYSTEM DESCRIBED 
BY AN EQUATION WITH CLOSED LINEAR MAPPING IN HILBERT SPACE}\\[5ex]
\textbf{Serhiy M.Zhuk}\\
Cybernetics faculty\\ Shevchenko National University of Kyiv, Ukraine\\beetle@unicyb.kiev.ua\\[15ex]
International Conference "Differential Equations and Topology", Moscow, 2008
\end{center}

\newpage
 In recent time the problem of the useful signal extraction from a noise
 becomes the question of the day. A simple example comes from the audio signal processing.
 When processing a piano signal $y(n)$ the sinusoidal components make it difficult to
detect the nonstationary mechanical noise $\eta(n)$ produced by the hammer
hitting the cords. One of possible ways to get rid of the noise $\eta(n)$ is to decompose the signal $y(n)$
modeled by
\begin{equation}
  \label{eq:yn}
y(n)=\varphi(n)+\eta(n)  
\end{equation}
into harmonic and noise parts. Mathematically this is performed by projecting
the signal $y(s)$ onto the signal and noise subspaces respectively.

In this example we use a simple additive uncertain model 
\begin{center}
  SIGNAL=USEFUL SIGNAL + NOISE
\end{center}
but in general case the signal may have more complicated structure. 
In this talk
we deal with the useful signal $\varphi$ modeled by the state of the linear
model
\begin{equation}
  \label{eq:Lfi}
  L\varphi=f
\end{equation}
while the measured signal $y$ is supposed to be a sum of the noise $\eta$
and a linear transformation of the useful signal $H\varphi$  
\begin{equation}
  \label{eq:Hfi}
y=H\varphi+\eta   
\end{equation}
We suppose here that $f$ is
unknown and belongs to the given set $G$. Here $G$ describes a type of
uncertainty. Another source of an uncertainty is represented by the
measurement noise $\eta$.  The problem of the useful
signal extraction is
\begin{center}\large
  to estimate a transformation $\ell(\varphi)$ of the state
$\varphi$  through an algorithm $\widehat{\ell(\varphi)}$ operating on
$y$. 
\end{center}

\newpage
 A state estimation framework for the linear dynamic models has several 
 unlike widely-used approaches: $H_2/H_\infty$ filtering, set-valued
  state observation and minimal guaranteed-cost (or minimax) estimation. $H_2$-estimators
 like Kalman or Wiener filters give the estimators of the signal generating
 process with minimum error variance also known as minimum variance
 filters. The optimal $H_\infty$ estimators minimize the norm of the operator
 that maps the unknown disturbances to the filtered errors. A set-valued
 estimators like Kurzhanskii's ellipsoidal informational sets observers 
describe a set of all possible realizations of the state $\varphi$ consistent with
 measurements $y$.

In this talk we address the
linear state estimation problem stated in terms of the minimax estimation
framework:
\begin{center}\large 
optimal linear algorithm $\widehat{\widehat{\ell(\varphi)}}$ minimizes the worst
case distance $
\sup_{\varphi,\eta}d(\ell(\varphi),\widehat{\ell(\varphi)})
$ \end{center}
In the sequel the
optimal algorithm $\widehat{\widehat{\ell(\varphi)}}$ is reffered to 
as a \emph{linear minimax estimation}.

If the noise $\eta$ is modeled by a random process then we apply a
\emph{linear minimax a priori estimation} approach. For 
the given estimation $\widehat{\ell(\varphi)}=(u,\cdot)+c$ we assign the worst
mean squared distance
\begin{equation}
  \label{eq:1}
  \sigma(\ell,u):=\sup_{L\varphi\in \mathscr{G},R_\eta\in\mathscr R}
    M(\ell(\varphi)-\widehat{\ell(\varphi)})^2
\end{equation}
and set
\begin{equation}
  \label{eq:2}
  \widehat{\widehat{\ell(\varphi)}}=(\hat u,y)+\hat c
\end{equation} where $(\hat u,\cdot)+\hat c$ has the minimal worst
mean squared distance
\begin{equation}
    \label{eq:huc}
    \hat\sigma(\ell):=\sigma(\ell,\hat u)=\inf_{u,c}\sigma(\ell,u), 
  \end{equation}
The minimal worst
mean squared distance $\hat\sigma(\ell)$ is called a \emph{minimax a priori
  error}. One can observe that the minimax a
priori estimation approach is focused on the minimax approximation of a linear
function $\ell(\varphi)$ on the convex set $L^{-1}(G)$ by means of an affine
function $u(y)$.

\newpage
If the noise $\eta$ has a deterministic nature we apply a minimax a posteriori
approach. Its key idea is a
geometrical one: the \emph{linear
  minimax a posteriori estimation in the direction $\ell$} is a map which takes
the perturbed information $y$ to the center of the interval
$\ell(G_y)$,
\begin{equation}
  \label{eq:Gy}
   G_y=\{\varphi:(L\varphi,y-H\varphi)\in\mathscr G\}
\end{equation} 
thus minimizing the worst case distance 
\begin{equation}
\rho(\varphi)=\sup_{\psi\in\mathscr G_y}
|(\ell,\varphi)-(\ell,\psi)|
\end{equation} 
The number 
\begin{equation}
\hat\rho(\ell):=\rho(\hat\varphi)=\inf_{\mathscr G_y}\rho(\varphi)
 \end{equation}
is called a \emph{minimax a posteriori error in the direction $\ell$.} 
 Observe, that a posteriori estimation approach is
 not applicable if the noise $\eta$ is modeled by a random process with unknown but bounded correlation operator.

Note that \textbf{a key feature of the presented state estimation approach may be
described as follows: we fix a class of linear operators $L$, $H$; given any
pair $L,H$ from that class we describe a class $\mathcal L$ of all solution
operators $\ell$ such that the minimax error is finite.} This
allows to extend  minimax estimation approach developed by Alexander Nakonechny to the general class of 
linear differential-algebraic models 
\begin{equation}
  \label{eq:xy}
  \begin{split}
&\dfrac d{dt}Fx(t)=C(t)x(t)+f(t),Fx(t_0)=f_0  ,\\
&y(t)=H(t)x(t)+\eta(t),t_0\le t\le T
  \end{split}
\end{equation}
with uncertain parameters $f_0,f(\cdot),\eta(\cdot)$.   

\newpage
Now let us show some results for that case when $L$ is closed linear
operator which maps the linear dense subset $\mathscr{D}(L)$ of Hilbert space
$\mathcal{H}$ into Hilbert space $\mathcal{F}$,
$H\in\mathscr{L}(\mathcal{H},\mathcal{Y})$. 
\par\noindent\textbf{Assumption 1} 
\begin{center}
the sets $R(L),H(N(L))$ are closed
\end{center}
\par\noindent\textbf{Assumption 2} 
\begin{center}
the set $R(T)=\{[Lx,Hx],x\in\mathscr D(L)\}$ is closed
\end{center}
\par\noindent\textbf{Theorem 1} Let
\begin{equation}
  \label{eq:t1}
  G=\{f\in\mathcal F:(Q_1f,f)\le 1\},
\eta\in\{\eta:M(Q_2\eta,\eta)\le 1\},
\end{equation} 
and suppose that Assumption 1 or 2 holds. Then the minimax a priori error
$\hat\sigma(\ell)$ is finite iff 
\begin{equation}
  \label{eq:t12}
  \ell\in R(L^*)+R(H^*)
\end{equation} In this
case $\hat\sigma(\ell)=(\ell,\hat p)^\frac12$ and the unique minimax a priori
estimation is given by $\widehat{\widehat{\ell(\varphi)}}=(\hat u,y)$ where $\hat 
u=Q_2H\hat p$, $\hat p$ is any solution of the equation 
\begin{equation}
  \label{eq:qeuler}
  \begin{split}
    &L^*\hat z=\ell-H^*Q_2H\hat p,\\
    &L\hat p=Q^{-1}_1\hat z
  \end{split}
\end{equation}

Theorem 1 is based on the general duality principle 
\par\noindent\textbf{Theorem 2}
Suppose that $\eta\in\{\eta:M(Q_2\eta,\eta)\le 1\}$, $\ell\in R(L^*)+R(H^*)$,
$G$ is a convex closed bounded subset of $\mathcal F$ and
\begin{equation}
  \label{eq:RLintG}
\mathrm{int}G\cap R(L)\ne\varnothing  
\end{equation}
Then the minimax a priori estimation problem $$
\sup_{L\varphi\in \mathscr{G},R_\eta\in\mathscr R}
    M(\ell(\varphi)-\widehat{\widehat{\ell(\varphi)}})^2\to\inf_{u,c}
$$
is equal to the optimal control problem $$
(Q_2u,u)+\inf\{c(\mathcal G,z)|L^*z=\ell-H^*u\}\to\min_u
$$

For instance Teorem 2 gives Kalman's duality theorem for linear ordinary differential
equations. Note that there exist a number of examples showing that
Condition~\eqref{eq:RLintG} is significant.

\newpage
Let's consider a posteriori estimations. 
\par\noindent\textbf{Theorem 3}
Let
\begin{equation}
  \label{eq:qf}
  G=\{(f,v):(Q_1f,f)+(Q_2\eta,\eta)\le 1\},
\end{equation} and suppose that Assumption 1 or 2 holds. Then the minimax a
posteriori error $\hat d(\ell)$ in the direction $\ell$ is finite iff
\begin{equation}
  \label{eq:4}
  \ell\in R(L^*)+R(H^*)
\end{equation} In this case 
\begin{equation}
  \label{eq:apos-err}
  \hat{\rho}(\ell)=(1-(Q_2y,y-H\hat\varphi))^\frac 12\hat\sigma(\ell)
\end{equation}
and minimax a posteriori estimation is given by
\begin{equation}
  \label{eq:est1}
  \widehat{\widehat{\ell(\varphi)}}=(\ell,\hat\varphi)=(\hat u,y)
\end{equation}
 where $\hat\varphi$ is
any solution of the equation 
\begin{equation}
  \label{eq:apoest:amap:thr}
  \begin{split}
    &L^*\hat q=H^*Q_2(y-H\hat\varphi),\\
    &L\hat\varphi=Q^{-1}_1\hat q
  \end{split}
\end{equation}
If $\hat d(\ell)<+\infty$ for any $\ell$ then 
\begin{equation}
  \label{eq:vest}
  \inf_{\varphi\in G_y}\sup_{x\in G_y}\|\varphi-x\|=
    \sup_{x\in G_y}\|\hat\varphi-x\|=(1-(Q_2y,y-H\hat\varphi))^\frac 12
    \max_{\|\ell\|=1}\hat{\sigma}(\ell)
\end{equation}
so that $\hat\varphi$ is the center of the a posteriori set $G_y$. 
\begin{flushleft}\large\bf
  Serhiy Zhuk, "Minimax state estimation for a dynamic system described 
by an equation with closed linear mapping in a Hilbert space,"\textit{Ukrainian
  Mathematical Journal}, In press
\end{flushleft}

\newpage
To apply the general theory to the minimax state estimation for~\eqref{eq:xy}
we need to translate~\eqref{eq:xy} in operator language. 
We say that $x(\cdot)$ satisfies~\eqref{eq:xy} if $t\mapsto Fx(t)$ is totally
continuous function, $Fx(t_0)=f_0$ and its derivative $t\mapsto\dfrac
d{dt}Fx(t)$ almost everywhere on $[t_0,T]$ is equal to $t\mapsto
C(t)x(t)+f(t)$. Denote by $\mathbb W_F$ the space of all $x(\cdot)\in\mathbb
L_2[t_0,T]$ such that $t\mapsto Fx(t)$ is totally
continuous function. Let $\mathcal D$ be the map which maps
$\mathbb W_F$ into $\mathbb L_2[t_0,T]\times\mathbb R^n$ by the rule $$
\mathcal Dx(t)=(\dfrac d{dt}Fx(t)-C(t)x(t),Fx(t_0)), x\in \mathbb W_F
$$ We prove that $\mathcal D$ is closed, its adjoint $\mathcal D^*$ is defined by the rule $$
\mathcal D^*(z_0,z)(t)=-\dfrac d{dt}F'z(t)-C'(t)z(t), (z,z_0)\in\mathbb W_{F'},
$$ where $\mathbb W_{F'}=\{(z,z_0),z_0=F^+Fz(t_0)+d,F'd=0,F'z\in\mathbb
W_2^1[t_0,T]\}$. 

\newpage
\par\noindent\textbf{Lemma 1} Let  $
 L' C(t)R'=\bigl(\begin{smallmatrix}
    C_1&&C_2\\C_3&&C_4
  \end{smallmatrix}\bigr)$, where 
\begin{equation}
    \label{eq:FLR}
    F=L\Lambda R, F^+=R'\Lambda^+L',LL'=E_m,RR'=E_n,\Lambda=\bigl(
\begin{smallmatrix}
  D^\frac 12&&0_{r,n-r}\\0_{m-r,r}&&0_{m-r,n-r} 
\end{smallmatrix}\bigr) 
\end{equation}
If $$\sup_{1>\varepsilon>-1}\|(\varepsilon^2E+C'_4C_4)^{-1}C'_2\|_{mod}<+\infty,
\eqno(a)$$ then $\mathcal D$ has a closed range. 

Now we can apply Theorems 1,2,3 to the minimax estimation of
the inner product $$
\ell(x)=\int_{t_0}^T(\ell,x)dt
$$ for the descriptor
model~\eqref{eq:xy}. 
\par\noindent\textbf{Proposition 1} 
Let \begin{equation} \label{setG2}
  \begin{split}
    &\mathscr
    G_2=\{R_\eta:\int_{t_0}^T\mathrm{tr}(Q_2(t)R_\eta(t,t))dt\le1\},R_\eta(t,s)=M\eta'(t)\eta(s)\\ 
    &\mathscr G_1=
    \{(f(\cdot),f_0):(Q_0f_0,f_0)+\int_{t_0}^T(Q_1(t)f(t),f(t))dt\le 1\}
  \end{split}
\end{equation}
Suppose that $H(N(\mathcal D))$ is closed set, where $Hx(t)=H(t)x(t)$. 
Under the condition $(*)$ we have  $$
  \hat{u}(t)=Q_2(t)H(t)p(t),\quad
  \hat\sigma^2(\hat u)=\int_{t_0}^T(\ell(t),p(t))dt,
  \eqno(b)
$$ where  \begin{equation}\label{eq:ddtFpFz:d}
    \begin{aligned}
      &\dfrac d{dt}F'z(t)=-C'(t)z(t)+H'(t)Q_2(t)H(t)p(t)-\ell(t),\\
      &\dfrac d{dt}Fp(t)=C(t)p(t)+Q^{-1}_1(t)z(t),\\
      &F'z(T)=0,Fp(t_0)=Q^{-1}_0(FF^+z(t_0)+d),F'd=0
  \end{aligned}
\end{equation}

\newpage
Lets consider a posteriori estimations. 
\par\noindent\textbf{Proposition 2} 
Let  $$
G=\{(f,\eta):\int_{t_0}^T(\|f(t)\|^2+\|\eta(t)\|^2)dt\le1\}\eqno(g)
$$ Then for any $\ell$ the minimax a posteriori estimation of $$
\ell(x)=\int_{t_0}^T(\ell,x)dt
$$ by observations $$
y(t)=x(t)+\eta(t),t_0\le
t\le T$$ is given by $$
\widehat{\widehat{\ell(x)}}=\int_{t_0}^T(\ell,\hat x)dt
\eqno(d)$$
where $$
\hat x(t)=[x_1(t),x_2(t)]
$$ and \begin{equation}\label{eq:x12q12:mnmx}
  \begin{split}
    &\dot x_1(t)=(C_1-C_2(E+C_4'C_4)^{-1}C'_4C_3)x_1(t)+\\
    &(C_2(E+C_4'C_4)^{-1}C'_2+E)q_1(t)+\\
    &C_2(E+C_4'C_4)^{-1}y_2(t), x_1(t_0)=0,\\
    &\dot q_1(t)=(-C'_1+C'_3C_4(E+C_4'C_4)^{-1}C'_2)q_1(t)+\\
    &C'_3C_4(E+C_4'C_4)^{-1}y_2(t)-y_1(t)+\\
    &(C'_3(E-C_4(E+C_4'C_4)^{-1}C'_4)C_3+E)x_1(t),q_1(T)=0,\\
    &x_2(t)=-(E+C_4'C_4)^{-1}C'_4C_3x_1(t)+\\
    &(E+C_4'C_4)^{-1}(C'_2q_1(t)+y_2(t)),\\
    &q_2(t)=-(E-C_4(E+C_4'C_4)^{-1}C'_4)C_3x_1(t)-\\
    &C_4(E+C_4'C_4)^{-1}(C'_2q_1(t)+y_2(t))
  \end{split}
\end{equation}
In this case we can estimate the whole state vector of~\eqref{eq:xy} with no
assumptions about structure of~\eqref{eq:xy} matrices. 
\begin{center}\large\bf
Serhiy Zhuk. "Properties of the linear mapping, induced by 
singular differential equations," \textit{Nonlinear oscillations}, No.~4, 
(2007), http://arxiv.org/abs/0705.3365
\end{center}

\newpage
Consider a linear time-variant system described by the following discrete-time
descriptor model
\begin{equation}
  \label{eq:state}
  \begin{split}
    &F_{k+1}x_{k+1}-C_k x_k = f_k, F_0x_0=q,\\
    &y_k=H_k x_k+g_k, k=0,1,\dots,
  \end{split}
\end{equation}
We assume that system input, measurement noise along with initial condition belong to the given set $$
 G=\{(q,\{f_k\},\{g_k\}):G(q,\{f_k\},\{g_k\})\leqslant1\}
$$ where 
\begin{equation*}
  \label{eq:Ginft}
  G(q,\{f_k\},\{g_k\})=(S q,q)+
  \sum_0^{\infty}(S_i f_i,f_i)+(R_i g_i,g_i),
\end{equation*}
\par\noindent\textbf{Theorem 4} 
 For the minimax a posteriori error
$\hat{\sigma}(\ell,N)$ in the direction $\ell$ to be finite it is necessary and sufficient to have 
  \begin{equation}
    \label{eq:erfn}
 Q_N^+Q_N\ell=\ell
  \end{equation}
Under this condition we have 
\begin{equation}\label{eq:err}
    \hat{\sigma}(\ell,N)=
    [1-\alpha_N+(Q^+_Nr_N,r_N)]^\frac 12(Q^+_N\ell,\ell)^\frac 12
  \end{equation}
 and 
   \begin{equation}
     \label{eq:est}
     \widehat{\widehat{(\ell,x_N)}}=(\ell,Q_N^+r_N)
   \end{equation}
where
   \begin{equation}
    \label{eq:Qk}
    \begin{split}
      &Q_k=
      H'_kR_kH_k+F'_k[S_{k-1}-S_{k-1}C_{k-1}P_{k-1}^+C'_{k-1}S_{k-1}]
      F_k,\\
      &Q_0=F'_0SF_0+H'_0R_0H_0,P_k=Q_{k}+C'_{k}S_{k}C_{k}
    \end{split}
  \end{equation}
$k\mapsto r_k$ denotes a recursive map that takes each natural number $k$ to the vector 
\begin{equation}
    \label{eq:rk1}
    \begin{split}
      &r_k=F'_kS_{k-1}C_{k-1}P_{k-1}^+
      r_{k-1}+H'_kR_k y_k,\\
      &r_0=H'_0R_0y_0
    \end{split}
\end{equation}
and 
\begin{equation}
  \label{eq:alphak}
  \begin{split}
    &\alpha_i=\alpha_{i-1}+(R_iy_i,y_i)-(P_{i-1}^+
 r_{i-1},r_{i-1}),\\
 &\alpha_0=(Sg,g)+(R_0y_0,y_0)
  \end{split}
\end{equation}
Observe that $$
(\ell,Q_N^+r_N)=0
\quad\text{if}\quad
(I-Q_N^+Q_N)\ell=\ell
$$ In this case
$$\ell(\mathscr G_y)=[-\infty,+\infty]
$$ and the estimation algorithm gives the
"center" of the real line or trivial estimation. Thus we can say that
\begin{center}\large\bf
the state $x_N$ is observable in the direction $\ell$ iff $\ell\in R(Q_N)$
\end{center} 
In other words we can assign nontrivial estimation
$\widehat{\widehat{(\ell,x_N)}}$ to the projection of the state  $x_N$ onto
the direction $\ell\in R(Q_N)$. 
\par\noindent\textbf{Definition 1}
The map $N\mapsto I_N=\mathrm{rank Q_N}$ is called \emph{an index of
  noncausality} of the model
\begin{equation*}
  \begin{split}
    &F_{k+1}x_{k+1}-C_k x_k = f_k, F_0x_0=q,\\
    &y_k=H_k x_k+g_k, k=0,1,\dots,
  \end{split}
\end{equation*}
 The index of noncausality $I_N$ gives the
number of linear independent observable directions in the system state
space. If it is equal to the state space dimension we say that the
system~\eqref{eq:state} is causal.

\newpage
Now we reveal the relation between the set-valued estimation approach and the
minimax a posteriori estimation method. Let $\pi_N$ be the map which takes each
vector $\mathbf x=(x_1\dots x_N)$ to its $N$-th component $x_N$, let
$\lambda_i(N)$ be the $i$-th eigenvalue of $Q_N$. Denote by $P_N$ the set-valued
map which takes each natural $N$ to the set $P_N$ of all $x\in\mathbb R^n$ such
that $x=\pi_N(\{x_k\})$, where $\{x_k\}_0^N$ satisfies 
\begin{equation}\label{eq:FkHk}
  \begin{split}
    &F_{k+1}x_{k+1}-C_k x_k = f_k, F_0x_0=q,\\
    &y_k=H_k x_k+g_k, k=0,1,\dots,
  \end{split}
\end{equation}
for some $q,\{f_k\}_0^{N-1}$, $\{g_k\}_0^N$. One can see that
$P_N$ is a vertical section of the model~\eqref{eq:state} trajectories set in
the system state space consistent with the measured signals $(y_1\dots y_N)$. Now we present an efficient
description of $P_N$ and its center for the causal descriptor model~\eqref{eq:FkHk}.
\par\noindent\textbf{Proposition 3}
 Suppose $I_N=n$. For any natural $N$ 
  \begin{equation}
  \label{eq:PnGy}
  P_N=\{x\in\mathbb R^n:(Q_Nx,x)-2(Q_N\hat x_N,x)+\alpha_N\leqslant1\}
\end{equation}
and
  \begin{equation}
    \label{eq:mobs}
    \begin{split}
      &\min_{x\in P_N}\max_{\tilde{x}\in P_N}\|x-\tilde{x}\|=\\
      &\max_{x\in P_N}\|x-\hat x_N\|=
      \frac{[1-\alpha_N+(Q_N\hat x_N,\hat x_N)]^\frac
        12}{\min_i\{\lambda^\frac 12_i(N)\}} 
    \end{split}
  \end{equation}
so that the central point of $P_N$ is given by $\hat{x}_N=Q_N^+r_N$. 
\begin{center}\large\bf
Zhuk S.M., "Recursive state estimation for noncausal
                       discrete-time
                       descriptor systems under uncertainties",
                       \textit{Automatica}, in press, 
                       http://arxiv.org/abs/0711.1334
\end{center}
\end{document}